\documentclass{amsart}
\usepackage{amsmath, amscd, amssymb, amsthm}
\usepackage{bbm}
\usepackage{latexsym}
\usepackage{amsfonts}
\usepackage{graphicx}
\usepackage[all,cmtip]{xy}
\usepackage[colorlinks,linkcolor=blue,breaklinks=blue,urlcolor=blue,citecolor=blue,anchorcolor=blue,pagebackref]{hyperref}%
\usepackage{geometry}
\newtheorem{theorem}{Theorem}

\newtheorem{question}[theorem]{Question}
\newtheorem*{example}{Example}

\newtheorem{conjecture}{Conjecture}

\renewcommand*\backref[1]{}
\renewcommand*\backrefalt[4]{ \ifcase #1 \or (cited on page #2) \else (cited on pages #2) \fi}

\newcommand{\be}{\begin{equation}}
\newcommand{\ee}{\end{equation}}
\newcommand{\bea}{\begin{eqnarray}}
\newcommand{\eea}{\end{eqnarray}}

\newcommand{\vs}{\vspace{0.5cm}}

\def\XXint#1#2#3{{\setbox0=\hbox{$#1{#2#3}{\int}$ }
\vcenter{\hbox{$#2#3$ }}\kern-.6\wd0}}

\begin{document}

\title[Constant holomorphic sectional curvature conjecture nd Fino-Vezzoni conjecture]{Constant holomorphic sectional curvature conjecture and Fino-Vezzoni conjecture}

\author{Fangyang Zheng}\thanks{Zheng is partially supported by National Natural Science Foundations of China
with the grant numbers 12141101 and 12471039.}
\address{Fangyang Zheng. School of Mathematical Sciences, Chongqing Normal University, Chongqing 401331, China}
\email{20190045@cqnu.edu.cn; franciszheng@yahoo.com} \thanks{}

\subjclass[2020]{53C55 (primary), 53C05 (secondary)}
\keywords{holomorphic sectional curvature, Hermitian space forms, Fino-Vezzoni conjecture, balanced metrics, pluriclosed metrics.}

\begin{abstract}
In this short essay, we will survey on two conjectures in non-K\"ahler geometry: the constant holomorphic sectional curvature conjecture and the Fino-Vezzoni conjecture. We aim at the broad audience and assume no expertise in non-K\"ahler geometry. We will discuss the history and recent developments on these two typical conjectures in the field.
\end{abstract}

\maketitle


\markleft{Fangyang Zheng}
\markright{Two conjectures in non-K\"ahler geometry}

\section{Introduction}\label{intro}

The topic of this survey lies in the field of {\em complex geometry}, which is an area of mathematics to study the  geometric, topological, and function theoretic  properties of complex manifolds, using calculus and complex analysis as primary tools.

{\em Complex manifolds} are high dimensional generalization of Riemann surfaces, which are $1$-dimensional complex manifolds (or complex curves). Given any complex manifold $M^n$ of complex dimension $n$, there always exist {\em Hermitian metrics} $g$ on it, which are Riemannian metrics compatible with the almost complex structure $J$ of $M^n$, namely, $g(Jx,Jy) =g(x,y)$ for any vector fields $x$ and $y$ on $M^n$. The pair $(M^n,g)$ is called a {\em Hermitian manifold}.  

Since the underlying topological space of $M^n$ is a $2n$-dimensional differentiable manifold, a Hermitian manifold $(M^n,g)$ is also a Riemannian manifold of real dimension $2n$. Denote by $\nabla$ the {\em Levi-Civita connection} (also called {\em Riemannian connection} sometimes), which is the unique linear connection on $M^n$ that is {\em metric} (namely, $\nabla g=0$) and {\em torsion-free} (namely, the associated torsion tensor $T(x,y):=\nabla_xy-\nabla_yx-[x,y]=0$ vanishes identically). 

For a Hermitian manifold $(M^n,g)$, if the Levi-Civita connection $\nabla$ is compatible with $J$ in the sense that $\nabla J=0$, then we say that $g$ is a {\em K\"ahler metric}, and call $(M^n,g)$ a {\em K\"ahler manifold}. If $\nabla J$ is not identically zero, then $g$ is said to be {\em non-K\"ahler}, and $(M^n,g)$ is called a {\em non-K\"ahler manifold}. For K\"ahler manifolds, Riemannian geometric studies fit perfectly with the complex analytic theory, and since the 1950s, K\"ahler geometry has seen great developments. 

An important class of compact complex manifolds are {\em projective manifolds}, namely, compact complex submanifolds of the complex projective space ${\mathbb C}{\mathbb P}^m$ for some large $m$. A famous theorem of Wei-Liang Chow says that such a manifold is always the common zero locus of a finite number of homogeneous polynomials of $m\!+\!1$ variables. Since the standard Hermitian metric on ${\mathbb C}{\mathbb P}^m$ (called the {\em Fubini-Study metric}) is K\"ahler, while by definition any complex submanifold of a K\"ahler manifold (when equipped with the restriction metric) is again K\"ahler, so all projective manifolds admit K\"ahler metrics. However, as we shall see later, there are lots of complex manifolds which do not admit any K\"ahler metrics, thus the study of general (non-K\"ahler) Hermitian geometry becomes necessary. In this survey, we will focus on two fundamental problems in non-K\"ahler geometry, the constant holomorphic sectional curvature conjecture and the Fino-Vezzoni conjecture.

Given any Hermitian manifold $(M^n,g)$, there are three canonical metric connections. Besides the aforementioned Levi-Civita connection $\nabla$, the {\em Chern connection} $\nabla^c$ is the unique metric connection such that $\nabla^cJ=0$ and its torsion $T^c$ has only $(2,0)$ part, while the {\em Bismut connection} $\nabla^b$ (also known as {\em Strominger connection} in some literature) is the unique metric connection such that $\nabla^bJ=0$ and $T^b$ is totally skew-symmetric, namely, $g(T^b(x,y),z) = - g(T^b(x,z),y)$ for any vector fields $x$, $y$, $z$ on $M$. The connection $\nabla^b$ was discovered independently by Bismut \cite{Bismut} and Strominger \cite{Strominger} in the 1980s, and we will call it Bismut connection from now on. When the metric $g$ is K\"ahler, these three canonical metric connections coincide: $\nabla =\nabla^c=\nabla^b$. When $g$ is not K\"ahler, then these three connections are mutually distinct, so we have three different kinds of geometry.  

\vspace{0.1cm}

Now let us discuss the constant holomorphic sectional curvature conjecture. First recall that in Riemannian geometry, the simplest kind of Riemannian manifolds are the so-called  {\em space forms,} which means complete Riemannian manifolds with constant sectional curvature. Their universal covers are respectively the sphere $S^n$, the Euclidean space ${\mathbb R}^n$, or the hyperbolic space ${\mathbb H}^n$, equipped with (scaling of) the standard metrics.

In the complex case, the sectional curvature of Hermitian manifolds in general can no longer be constant (unless it is identically zero, i.e., {\em flat}). Instead one requires the {\em holomorphic sectional curvature} to be constant. When the metric is K\"ahler, one gets the so-called {\em complex space forms,} namely complete K\"ahler manifolds with constant holomorphic sectional curvature. Analogous to the Riemannian case, their universal covers are the complex projective space ${\mathbb C}{\mathbb P}^n$, the complex Euclidean space ${\mathbb C}^n$, or the complex hyperbolic space ${\mathbb C}{\mathbb H}^n$, all equipped with (scaling of) the standard metrics.

When a Hermitian metric is not K\"ahler, its curvature tensor does not obey all the K\"ahler symmetries in general. As a result, the holomorphic sectional curvature could no longer determine the entire curvature tensor. So one would naturally wonder about when can the holomorphic sectional curvature be  constant. In this direction, a long-standing conjecture is the following:
\begin{conjecture}[Constant Holomorphic Sectional Curvature Conjecture] \label{conj1}
Given any compact Hermitian manifold, if the holomorphic sectional curvature of its Chern (or Levi-Civita) connection is a constant $c$, then when $c\neq 0$ the metric must be K\"ahler (hence it is a complex space form), while when $c=0$ the metric must be Chern (or Levi-Civita) flat. 
\end{conjecture}
 Note that when $n\geq 3$, there are  compact Chern flat (or Levi-Civita flat) manifolds that are non-K\"ahler. Compact Chern flat manifolds are compact quotients of complex Lie groups by the classic theorem of Boothby \cite{Boothby}, and compact Levi-Civita flat Hermitian threefolds were classified in \cite{KYZ}, while for dimension $4$ or higher there are lots of such metrics, yet the full classification is still unknown. 

The compactness assumption in the conjecture is a necessary one, without which the conclusion fails. For instance, in \cite{CCN}, Chen, Chen, and Nie constructed an example of a complete  Hermitian manifold which has vanishing Chern holomorphic sectional curvature, yet it is not Chern flat.

For $n=2$, Conjecture \ref{conj1} is true by the combined effort of Balas and Gauduchon \cite{Balas, BG}, Sato and Sekigawa \cite{SS}, and Apostolov, Davidov, and Muskarov \cite{ADM} in the 1980s and 1990s. In higher dimensions, the first substantial result towards this conjecture is the one obtained by Davidov, Grantcharov, and Muskarov \cite{DGM}, in which they showed among other things that the only twistor space with constant holomorphic sectional curvature is the complex space form ${\mathbb C}{\mathbb P}^3$. In \S 3, we will discuss some other partial results to Conjecture \ref{conj1}.

\vspace{0.1cm}

Next let us discuss the Fino-Vezzoni conjecture, which is about the non-coexistence of two special types of Hermitian structures. Let $\omega$ be the global $(1,1)$-form associated to the metric $g$ of a given  Hermitian manifold $(M^n,g)$, defined by $\omega (x,y) = g(Jx,y)$. As it is well-known, $g$ is K\"ahler if and only if $d\omega =0$. Since any complex manifold admits a Hermitian metric, the study of non-K\"ahler geometry are often focused on analyzing special types of Hermitian manifolds. Two types of special Hermitian structures that are widely studied are {\em balanced metrics} and {\em pluriclosed metrics}. The former is defined by $d(\omega^{n-1})=0$ where $n$ is the complex dimension of the manifold, while the latter is defined by $\partial \overline{\partial}\omega =0$. Both are obvious generalizations to K\"ahlerness. Balanced metrics were first introduced by Michelsohn \cite{Michelsohn} in 1982 where she gave a Harvey-Lawson type characterization for such metrics.  Pluriclosed metrics are also called {\em strong K\"ahler with torsion} metrics, abbreviated as SKT metrics. In \cite{FinoTomassini}, Fino and Tomassini gave an excellent survey  on this type of special metrics.

The balanced structure and pluriclosed structure seem to be mutually exclusive, and in \cite{FV15} and \cite{FV16}, Fino and Vezzoni proposed the following interesting conjecture:

\begin{conjecture}[Fino-Vezzoni Conjecture] \label{conj2}
If a compact complex manifold $M^n$ admits a balanced metric and a pluriclosed metric, then it must admit a  K\"ahler metric. 
\end{conjecture}

We remark that since the balanced condition coincides with the K\"ahler condition when $n=2$, Conjecture \ref{conj2} starts in complex dimension $3$. Also, if a Hermitian metric $g$ is simultaneously balanced and pluriclosed, then it is K\"ahler by a theorem of Alexsandrov and Ivanov \cite{AI}. So the point of the conjecture is that the two metrics involved could be different. The conjecture says that for non-K\"ahlerian compact complex manifolds, balanced ones and pluriclosed ones form disjoint subsets.  The conjecture has been drawing a lot of attentions in complex geometry since its inception, and so far it has been confirmed for many special types of Hermitian manifolds, yet even in complex dimension $3$ it is still open in its full generality. In \S 4 we will discuss these partial results towards this interesting conjecture. 

We refer the readers to the excellent survey papers by D. Angella \cite{Angella}, by A. Fino \cite{Fino}, and by A. Fino and G. Grantcharov \cite{FG}. In these survets the authors discussed a number of important issues in non-K\"ahler geometry, including the two conjectures above. This article is organized as follows. In \S 2, we will discuss with some generality the distinction between K\"ahler and non-K\"ahler manifolds, and also present several examples of non-K\"ahler manifolds that are well-known and of historical significance. In the next two sections, we will discuss some old and new partial results towards the two aforementioned conjectures.

\vspace{0.3cm}

\section{Examples of non-K\"ahler manifolds}

\subsection{Almost  complex manifolds}
Let us start with almost complex manifolds.  In the 1940s, Heinz Hopf and Charles Ehresmann realized that on any complex manifold, the complex structure naturally induces an endomorphism $J$ on the tangent bundle of the underlying differentiable manifold, satisfying the property that $J^2=-I$, namely, the composition of $J$ with itself equals to minus of the identity map. They called this $J$ an {\em almost complex structure}, and they call any $2n$-dimensional  differentiable manifold equipped with an almost complex structure an {\em almost complex manifold}. It is not difficult to see that the presence of $J$ implies that the manifold is orientable, so a naturally question becomes:

\vspace{0.05cm}

{\em Given an orientable $2n$-dimensional differentiable manifold $N^{2n}$, when does it admit an almost complex structure?}

\vspace{0.05cm}

The above question can be rephrased in terms of the structure group of the tangent bundle $TN$ of $N$, namely, when can the structure group of  $TN$ be reduced from $GL(2n,{\mathbb R})$ to  $GL(n,{\mathbb C})$? This is purely a question in algebraic topology,  and in their terminology, the question becomes when can the classifying map $t_N: N \rightarrow BO(2n)$ be lifted to  $\tilde{t}_N: N \rightarrow BU(n)$? The answer turns out to be the following: 

\vspace{0.05cm}

{\em An orientable differentiable manifold $N^{2n}$ admits an almost complex structure if and only if the obstruction classes  $\alpha_i=0$ for all $0\leq i\leq 2n-1$, where $\alpha_i \in H^{i+1}(N, \pi_i (Y_n))$ and $Y_n$ is the compact Hermitian symmetric space $Y_n=SO(2n)/U(n)$.}

\vspace{0.05cm}

Note that $Y_n$ is the set of all (constant) almost complex structures on the vector space ${\mathbb R}^{2n}$ compatible with the standard metric and orientation. We remark that, although the above statement gives a necessary and sufficient criterion to detect whether a given $N^{2n}$ is almost complex or not, in reality it is often not easy to use as these obstruction classes might be difficult to compute in general. This is why it is still an active area of current research to detect whether a given class of manifolds are almost complex or not. For instance, in a relatively recent work \cite{GK} in 2019, the authors showed that the $k \#{\mathbb C}{\mathbb P}^{2n}$ admits almost complex structure if and only if $k$ is odd. Here the symbol means the connected sum of $k$ copies of the complex projective space of complex dimension $2n$. For $n=1$, this is due to M. Audin in 1991, while for $n=2$, it is due to S. Muller and H. Geiges in 2000. 

On the positive side, it is well known that for even-dimensional sphere $S^{2n}$, the only almost complex ones are $S^2$ and $S^6$, while for product of two even-dimensional spheres, Sutherland showed that the only new addition is $S^2\!\times \!S^4$. Also, the above obstruction-theoretic criterion is very effective in low dimensional cases, especial in real dimensions $4$ and $6$. For real dimension $4$, we have the following famous result of Wen-Tsun Wu:

\begin{theorem}[Wu] 
A closed orientable $4$-manifold $N^4$ admits an almost complex structure $J$ if and only if there exists a class $h\in H^2(N, {\mathbb Z})$ such that
$ h^2=3\sigma +2\chi$ and $h\equiv w_2$ mod $2$, where $\sigma$ is the signature, $\chi$ is the Euler number, and $w_2$ is the second Stiefel-Whitney class. Furthermore, $h=c_1(N,J)$.
\end{theorem}

This result provides an effective way to check whether a given $4$-manifold is almost complex or not. Using this, along with the classification theory for compact complex surfaces, Van de Ven constructed in 1966 the first example of a closed $4$-manifold which is almost complex but not a complex manifold. Later, more such examples were discovered. As an illustration to the usage of Wu's Theorem, let us consider the following:

\begin{example} [{\bf 1}]
Consider the $4$-manifold  $N^4=(S^2\!\times\!S^2)\# (S^3\!\times\!S^1)\# (S^3\!\times\!S^1)$. It is an almost complex manifold but is not a complex manifold.
\end{example}

To see this, we note that signature, Euler number, and second Stiefel-Whitney class of $N$ are all $0$, so by taking $h=0$ in Wu's Theorem we know that $N^4$ is almost complex. On the other hand, it is not hard to see that the fundamental group $\pi_1(N)=F_2$, the free group of two generators. As it is well-known, $F_3$ is isomorphic to a subgroup of index $2$ in $F_2$, which means that $N$ admits an unbranched double cover $N'$ whose fundamental group is $F_3$. In particular, their first Betti numbers are respectively $b_1(N)=2$ and $b_1(N')=3$, as the first homology group is the abelianization of the fundamental group. Now we conclude that $N$ cannot be a complex manifold as we have the following well-known  result in complex surface theory:

\vspace{0.05cm}

{\em A compact complex surface is K\"ahlerian (namely it admits a K\"ahler metric) if and only if its first Betti number is even.}

\vspace{0.05cm}

If we assume that $N$ is a complex surface, then it would admit a K\"ahler metric as $b_1(N)$ is even. Then the metric can be lifted onto the double cover so $N'$ became a K\"ahler surface, which would contradict to the fact that $b_1(N')$ is odd.  Therefore we know that $N$ cannot be a complex manifold.  

\vspace{0.1cm}

In real dimension $6$, the aforementioned obstruction-theoretic criteria is still quite useful, as $Y_3={\mathbb C}{\mathbb P}^3$ and $\pi_i(Y_3)=0$ for all $i\leq 5$ except $\pi_2(Y_3)\cong {\mathbb Z}$. Furthermore, the obstruction class $\alpha_2\in H^3(N,\pi_2(Y_3))$ lies in the image of the Bockstein homomorphism $H^2(N,{\mathbb Z}_2) \rightarrow H^3(N,{\mathbb Z})$, so one has the following:

\begin{theorem}
Given a closed  orientable $6$-manifold $N$, if $H^3(N,{\mathbb Z})$ contain no $2$-torsion element, then $N$ is almost complex.
\end{theorem} 

D. Sullivan showed that any compact $6$-dimensional real hyperbolic space admits a finite unbranched cover whose $H^3$ does not have $2$-torsions, so the cover is always an  almost complex manifold.

\vspace{0.1cm}

In contrast to the $4$-dimensional situation, in higher dimensions there is no known example of closed $2n$-manifold which is almost complex but not complex, for any $n\geq 3$. In particular, a famous open question is whether $S^6$ admits a complex structure or not. Many people believe that it should not admit a complex structure, but there are also people who believe it does. In fact, Yau proposed the following conjecture in 1990s:

\begin{conjecture}[Yau]
For $n\geq 3$, any closed $2n$-manifold admitting an almost complex structure will also admit a complex structure. 
\end{conjecture}

When an almost complex structure $J$ is induced from a complex structure, then we say that $J$ is {\em integrable}. A necessary condition for $J$ to be integrable is that its Nijenhuis tensor $N_J$ must vanish identically. $N_J$ is defined by
$$ N_J(x,y) = [x,y] -[Jx,Jy] + J[Jx,y] + J[x,Jy], $$
where $x$ and $y$ are vector fields  on the manifold. Denote by $T^{1,0}$ the bundle of complex tangent vector fields of $(1,0)$ type, namely those in the form $x-\sqrt{-1}Jx$ where $x$ is real. Then clearly $N_J=0$ is equivalent to $[T^{1,0}, T^{1,0}] \subset T^{1,0}$. Thus when $J$ is induced form a complex structure, the Nijenhuis tensor must vanish. The converse to this is also true, which is a deep theorem of Newlander and Nirenberg in 1957:

\begin{theorem}[Newlander-Nirenberg]
An almost complex structure $J$ is integrable if and only if its Nijenhuis tensor $N_J=0$ identically. 
\end{theorem} 

We remark that on an almost complex manifold $(N,J)$, the set of all almost complex structures on $N$ forms an infinite dimensional space ${\mathcal A}{\mathcal C}(N)$, so to detect on whether or not it contains  an integrable element would be a daunting task. This  perhaps  explains why the conjecture about $S^6$ is still open today. 

\vspace{0.1cm}

\subsection{Compact K\"ahler manifolds}  First we remark that in complex geometry, people are primarily concerned with compact complex manifolds, as the non-compact ones are simply too complicated to deal with, except in some very special circumstances such as domains in ${\mathbb C}^n$, which are the main objects of study in Several Complex Variables. Non-compact manifolds not only can have infinite topology (e.g., their Betti numbers can be infinity), and even when their topology is simple, the complex structure can be very strange. For instance, in the 1950s E. Calabi constructed a complex structure on ${\mathbb R}^6$, and the resulted complex threefold $M^3$ does not have any non-constant global holomorphic function on it, which is in sharp contrast to ${\mathbb C}^3$. 

Compact complex manifolds can be divided into two disjoint subsets:  K\"ahler manifolds and non-K\"ahler manifolds. The former means those which admit a K\"ahler metric (strictly speaking they should be called {\em K\"ahlerian manifolds} but often time people use the shorter name), and the latter means those which do not admit any K\"ahler metric.  K\"ahler manifolds have been extensively studied since the 1950s, under the leadership of  great mathematicians such as Chern, Hodge, Serre, Kodaira, Atiyah,  Hirzebruch, Calabi, H\"ormander,  Yau, Siu, Mori, Donaldson, Demailly, and others. The study of non-K\"ahler manifolds, on the other hand, have been lagged behind until late last century and the new millennium. The reason for this is perhaps two-fold. On one hand, as we mentioned before, K\"ahler manifolds are much easier to analyze at least from the differential geometric point of view, where the Levi-Civita connection respect the almost complex structure so Riemannian geometry is in sync with the complex analysis. On the other hand, the set of compact K\"ahler manifolds contains all projective manifolds, which is the main concern of algebraic geometers.  

Note that the set of all compact K\"ahler manifolds is somewhat close to its subset of all projective manifolds. For instance, the following is still an open conjecture in complex geometry:

\begin{conjecture}[Peternell] \label{conj4}
Any minimal K\"ahler variety has an algebraic approximation.
\end{conjecture}

The conjecture implies in particular that any minimal compact K\"ahler manifold can always  be deformed (namely, the complex structure can be  continuously changed) into a projective manifold. Here the word `minimal' means that the manifold cannot be blown down any further. 

For a long time people believed that the statement should be true even without the `minimal' assumption, namely, {\em any compact K\"ahler manifold can be deformed into a projective one}. This last statement is called the {\em Kodaira problem} in the literature. Kodaira proved that it is true in complex dimension $2$, namely, any compact K\"ahler surface can always be deformed to a projective surface. However, in 2004, Clair Voison constructed counterexamples to the statement in dimensions $n\geq 4$  (\cite{Voisin04}, \cite{Voisin06}). 

At this point the readers would naturally be curious about the $n=3$ situation. A recent preprint by H-Y Lin \cite{Lin17} gives an affirmative answer to the Kodaira problem for $n=3$.  

The Hodge and Lefschetz decomposition theorems imply that compact K\"ahler manifolds obey some very strong topological restrictions. For instance, their Betti numbers satisfy  $b_{2k}>0$ and $b_{2k-1}$ is even, for any $1\leq k\leq n$ where $n$ is the complex dimension of the manifold. Also, a famous result by Deligne, Griffiths, Morgan, and Sullivan in the 1970s says that the rational homotopy type of a compact K\"ahler manifold is completely determined by its cohomology ring. 

Of course there is still a big distinction between projective manifolds and compact K\"ahler manifolds in general. The famous Hodge Conjecture, which states that any rational $(p,p)$-class $\alpha \in H^{p,p}(M)\cap H^{2p}(M, {\mathbb Q})$ of a projective manifold $M^n$ can be represented by an analytic cycle, for any $1\leq p\leq n$. Clair Voisin showed in her stunning work in 2002 that the conjecture fails for compact K\"ahler manifolds. 

Nonetheless, compact K\"ahler manifolds are still `close' to projective manifolds, and people often try to generalize results that are valid for projective manifolds to the K\"ahler case. An important example is the  result of Demailly and Paun in 2004 on the numerical characterization of the K\"ahler cone of a compact K\"ahler manifold, which is the K\"ahler version of the classic Nakai-Moishezon criterion for projective manifolds. 

Another example is the study of {\em K\"ahler groups}, which means the set of isomorphism classes of fundamental groups of compact K\"ahler manifolds. Denote by ${\mathcal K}_n$ (or ${\mathcal P}_n$) the set of all isomorphism classes of fundamental groups of compact K\"ahler manifolds (or projective manifolds) of complex dimension $n$. Then clearly ${\mathcal P}_n \subset {\mathcal K}_n$ for each $n$. On the other hand, since $\pi_1(M)\cong \pi_1(M\!\times\!{\mathbb C}{\mathbb P}^1)$, we have ${\mathcal P}_2 \subset {\mathcal P}_3 \subset {\mathcal P}_4\subset \cdots $ and ${\mathcal K}_2 \subset {\mathcal K}_3 \subset {\mathcal K}_4 \subset \cdots $. 

The union of all ${\mathcal K}_n$, namely the set
$ {\mathcal K} := \bigcup_{n=2}^{\infty}  {\mathcal K}_n $,
is called the K\"ahler group. Given any projective manifold $M^n$,  let $D\subset M$ be a smooth ample divisor. Then the Lefschetz hyperplane section theorem says that $\pi_1(D)\cong \pi_1(M)$ if $n\geq 3$. Therefore we know that ${\mathcal P}_2 = {\mathcal P}_n$ for any $n\geq 2$, and this set will be denoted by  $ {\mathcal P}$ and called the {\em projective group}. By Kodaira's theorem, any compact K\"ahler surface can be deformed to a projective one, in particular they are diffeomorphic to each other, thus ${\mathcal K}_2={\mathcal P}_2$. For $n=3$, Claudon, H\"oring, and Lin showed in 2019 \cite{CHL} that ${\mathcal K}_3={\mathcal P}_3$. They even conjectured that the same should hold in higher dimensions, namely, ${\mathcal K}={\mathcal P}$. This shows that, at least in complex dimension $3$,  compact K\"ahler manifolds and projective manifolds share the same set of fundamental groups. From this point of view, there are a lot more non-K\"ahler ones than K\"ahler ones, as we shall see by the twistor spaces example in the next subsection. 

\vspace{0.1cm}

\subsection{Non-K\"ahler manifolds} In this subsection we will recall some well-known constructions of non-K\"ahler complex structures. As we mentioned before, a complex manifold $M^n$ of complex dimension $n$ can be viewed as an (oriented) differentiable manifold $N^{2n}$ of real dimension $2n$ equipped with a complex structure, which means an integrable almost complex structure. If we fix the underlying smooth manifold $N^{2n}$ and consider the set of all almost complex structures on $N$:
$$ {\mathcal A}{\mathcal C}(N) = \{ J \in \mbox{End}(TN) \mid J^2=-I \}. $$
Write ${\mathcal C}(N) \subset {\mathcal A}{\mathcal C}(N) $ for the subset of all $J$ that are integrable, namely with $N_J=0$. Both  ${\mathcal A}{\mathcal C}(N) $ and  ${\mathcal C}(N)$ can be viewed as manifold of infinite dimension. If $f$ is a diffeomorphism from $N$ onto itself, then one can define  the pull back  $f^{\ast}(J):=(f^{-1})_{\ast}\circ J \circ f_{\ast}$ of an almost complex structure $J$, so the (infinite dimensional) diffeomorphism group $\mbox{Diff}(N)$ acts on both ${\mathcal A}{\mathcal C}(N) $ and  ${\mathcal C}(N)$. In somewhat loose terms, one may call the quotients 
$$ {\mathcal M}_a(N)= {\mathcal A}{\mathcal C}(N)/\mbox{Diff}(N) \ \ \ \ \ \mbox{and} \ \ \ \ \ {\mathcal M}(N)= {\mathcal C}(N)/\mbox{Diff}(N) $$
the {\em moduli space of almost complex structures} on $N$ or the {\em moduli space of complex structures} on $N$, respectively. Because of the complexity of the diffeomorphism group and its action, it is in general very difficult to analyze $ {\mathcal M}_a(N)$ or $ {\mathcal M}(N)$ this way, despite its natural appealing.  

Kodaira and Spenser developed the {\em deformation theory} for complex structures. Let $\pi : X \rightarrow B$ be a  holomorphic submersion (namely, surjective and without any critical point) between complex manifolds, where $B\subset {\mathbb C}^m$ is a small ball, and $M^n=\pi^{-1}(0)$ is a compact complex manifold. By Ehresmann's fibration theorem, $X$ is diffeomorphic (but in general not biholomorphic) to $M\times B$, in other words,  $M_t=\pi^{-1}(t)$  is diffeomorphic (but may not be biholomorphic) to $M$, for any $t\in B$. If we denote  the underlying smooth manifold of all $M_t$ by $N$, and denote the (integrable almost) complex structure of $M_t$ by $J_t$, then $\{ J_t \}_{ t\in B}$ is a family of deformation of the complex structure $J_0$ on the central fiber $M$. An important theorem of Kodaira states that if $M$ is K\"ahler, then $M_t$ will be K\"ahler for all $t$ with $|t|$ sufficiently small. This is known as the stability of K\"ahlerness.  

A surjective holomorphic map $f: M\rightarrow N$ between two compact complex manifolds of the same dimension is called a {\em modification} if there exists a proper analytic subset $Z\subset N$ such that the restriction of $f$ on $M\setminus f^{-1}(Z) \rightarrow N\setminus Z$ is a biholomorphism. Any bimeromorphic map  $f: X\rightarrow Y$ can be written as the composition $f=g\circ h^{-1}$ where both $g:Z\rightarrow Y$ and $h: Z\rightarrow X$ are modifications. Let $Z\subset X$ be a compact complex submanifold of codimension at least $2$. Denote by $\widetilde{X}$ the blowing-up of $X$ along the center $Z$. Then the projection map $\pi : \widetilde{X} \rightarrow X$ is called the blowing-down, and its inverse (which is only a meromorphic map) is called the blowing-up. Blowing-downs are very special examples of modifications. A fundamental theorem of Hironaka (resolution of singularities) states that any bimeromorphic map can be written as the composition of finitely many blowing-ups and blowing-downs, all along smooth centers.  

Note that bimeromorphic manifolds have the same fundamental group hence the same first Betti number, so for compact complex surfaces, K\"ahlerness is preserved under deformation or bimeromorphic (birational) changes. However, in dimensions $3$ or higher, K\"ahlerness is not preserved under either deformation or bimeromorphism. 

\begin{example}[{\bf 2}] In his PhD thesis in 1960, Hironaka constructed examples of a family of compact complex  threefolds $X^3_t$, parametrized by $t\in [0,1]$. For all $t\in [0,1)$, $X_t$ is projective, but $X_1$ is not K\"ahler. This shows that K\"ahlerness (or projectiveness) is not a closed condition under deformation. He also showed that a blowing up of $X_1$ is projective, so the blowing-down of a K\"ahler manifold may no longer be K\"ahler.
\end{example}

Note that the blowing-up (along smooth centers) of a K\"ahler manifold is always K\"ahler, as one can add the pull-back of the  K\"ahler form downstair by a suitable $\partial \overline{\partial}$-exact term which is supported in a neighborhood of the exceptional divisor so that the sum is strictly positive.  Similarly, the blowing-up of a projective manifold is always projective. 

Hironaka's example shows that compact K\"ahler manifolds can be deformed into a non-K\"ahler one, and the K\"ahlerness  condition is not closed under bimeromorphic (birational) equivalence.  

\begin{example}[{\bf 3}] A (primary) {\bf \em Hopf manifold} is a quotient of  ${\mathbb C}^n\setminus \{ 0\}$ (where $n\geq 2$) by the infinite cyclic group ${\mathbb Z}$ generated by a holomorphic contraction. The complex structure on the punctured complex Euclidean space descends down onto the quotients, which are all diffeomorphic to $S^1\times S^{2n-1}$. The latter has $b_1=1$ which is odd, so all Hopf manifolds are non-K\"ahler.
\end{example}

For $n=2$, Kodaira classified Hopf surfaces in \cite{Kodaira}, where the primary ones are either linear (Class $1$) $M^2_{a,b}=({\mathbb C}^2\setminus \{ 0\})/{\mathbb Z}f_{a,b}$, or non-linear (Class $0$) $M^2_{a, \lambda , m}=({\mathbb C}^2\setminus \{ 0\})/{\mathbb Z}h_{a, \lambda , m}$, with the corresponding holomorphic contraction given respectively by
$$ \left\{  \begin{split} f_{a,b}(z_1, z_2)=(a z_1, b z_2), \ \ \ \ \ \ \ \ a,b\in {\mathbb C}, \ 0<|a|\leq |b|<1; \hspace{4.25cm} \\ h_{a,\lambda, m}(z_1, z_2)= (a^mz_1, \,az_2+\lambda z_1^m), \ \ \ \  a, \lambda \in {\mathbb C}, \  m\in {\mathbb N}, \ 0<|a|<1, \  \lambda \neq 0, \ m\geq 2. \end{split} \right. $$


\begin{example}[{\bf 4}] In 1953, E. Calabi and B. Eckmann \cite{CalabiEckmann} constructed a family of complex structures on the product $S^{2p+1}\!\times S^{2q+1}$ of two odd dimensional spheres. Here $p$, $q$ are positive integers. All Calabi-Eckmann manifolds are non-K\"ahler, since its second Betti number $b_2=0$.
\end{example}

Let  $p$, $q$ be positive integers, and fix a constant $t\in {\mathbb C}\setminus {\mathbb R}$. Consider the holomorphic action of ${\mathbb C}$ on $({\mathbb C}^{p+1}\setminus \{0\}) \!\times \! ({\mathbb C}^{q+1}\setminus \{0\}) $ given by $z(x,y)= (e^zx, e^{tz}y)$. The action is free and proper, so the quotient $M$ becomes a complex manifold of complex dimension $n=p\!+\!q\!+\!1$, which is diffeomorphic to $S^{2p+1}\!\times \!S^{2q+1}$.

Note that if we denote by $f_p: S^{2p+1} \rightarrow {\mathbb C}{\mathbb P}^p$ the Hopf fibration map, which is the restriction of the  natural projection $\pi_p :{\mathbb C}^{p+1}\setminus \{0\} \rightarrow {\mathbb C}{\mathbb P}^p$ on the unit sphere. Then $f_p\!\times\!f_q$ makes $M$ a holomorphic fiber bundle over ${\mathbb C}{\mathbb P}^p\!\times \! {\mathbb C}{\mathbb P}^q$, with fiber being the elliptic curve $E_t={\mathbb C}/({\mathbb Z}+t{\mathbb Z})$. It is actually an $E_t$-principal bundle.


\begin{example}[{\bf 5}] In 1953, H. Samelson \cite{Samelson} constructed complex structures on any even-dimensional compact Lie groups. They are now known as the {\bf \em Samelson spaces}. More generally, H-C Wang \cite{Wang} gave a systematic study on compact complex homogeneous spaces. 
\end{example}


An important class of compact complex threefolds are given by the twistor spaces of Roger Penrose.

\begin{example}[{\bf 6}] Let $(X^4,g)$ be a closed orientable Riemannian $4$-manifold. Let $\pi : M\rightarrow X$ be the $S^2$-bundle where the fiber $\pi^{-1}(x)$ at $x\in X$ is the set of all almost complex structures on $T_xX\cong {\mathbb R}^4$ compatible with the metric and orientation.  Given any $y=(x,v)\in M$, where $x\in X$ and $v\in \pi^{-1}(x)\cong S^2$, the tangent space of $M$ splits as the direct sum of the horizontal space $H_y\cong T_xX$ and the vertical space $V_y\cong T_vS^2$. The fiber $S^2$ is the rational curve, hence has a complex structure $J_0$, and let $J_v$ be the almost complex structure on $H_y$ corresponding to $v\in S^2$. Then $J=J_v \!+\! J_0$  naturally becomes an almost complex structure  on $M$. There is also a (conformal family of) naturally induced  metric $h_t$ on $M$, where $t>0$ is any positive real number, so that $(M,J,h_t)$ becomes an almost Hermitian manifold. 

In 1978, Atiyah, Hitchin, and Singer \cite{AHS} showed that $J$ will be integrable if and only if $g$ is anti-self-dual. In this case $M$ becomes a compact complex threefold, called the {\bf \em twistor space} of $(X,g)$. 
\end{example}

By the famous gluing theorem of C. Taubes, given any closed $4$-manifold $X$, denote by $\widetilde{X}=X\# k\overline{{\mathbb C}{\mathbb P}^2} $ the connected sum of $X$ with $k$ copies of  ${\mathbb C}{\mathbb P}^2$ (with reversed orientation), then $\widetilde{X}$ always admits an anti-self-dual metric when $k$ is sufficiently large.  Since any finitely presented group can be the fundamental group of a closed oriented $4$-manifold, we conclude that

\vspace{0.1cm}

{\em Any finitely presented group can be the fundamental group of a compact complex threefold.}

\vspace{0.1cm}

This is in sharp contrast to the K\"ahler case, as the set ${\mathcal K}_3$ of fundamental groups of compact K\"ahler threefolds, is equal to the set  of projective groups ${\mathcal P}$. 

The twistor spaces form an interesting and important class of compact complex manifolds. They are often used as testing ground for questions and conjectures in non-K\"ahler geometry. The following results are known for them:

\begin{itemize}
\item For any $t>0$ the metric $h_t$ is balanced on $(M,J)$. 
\item (Hitchin, \cite{Hitchin81}) All twistor spaces are non-K\"ahlerian except two of them: ${\mathbb C}{\mathbb P}^3$ which is the twistor space of the standard $S^4$, and the flag threefold ${\mathbb P}(T_{{\mathbb C}{\mathbb P}^2} )$, which is the twistor space of the standard ${\mathbb C}{\mathbb P}^2$. 
    \item (Verbitsky \cite{Verbitsky}) All twistor spaces are rationally connected, meaning that given any two points on the threefold there is a rational curve passing through both of them.  
\end{itemize}

\vspace{0.3cm}

\section{Constant holomorphic sectional curvature conjecture}

In this section, we will discuss the developments surrounding the constant holomorphic sectional curvature conjecture stated as Conjecture \ref{conj1} in the introductions. In shorts, the conjecture seeks to understand the Hermitian version of (compact) spaces forms, namely, compact Hermitian manifolds with constant Chern (or Riemannian) holomorphic sectional curvature. Denote this constant by $c$. The conjecture states that, when $c\neq 0$, then the metric must be K\"ahler thus the manifold is a complex space form, while when $c=0$, then the metric must be Chern (or Riemannian) flat. Such a metric can be non-K\"ahler when the dimension is $3$ or higher.

As mentioned in the introduction section, Conjecture \ref{conj1} is true in the $2$-dimensional case. This is due to a combined effort by several groups of authors.  For the $c<0$ and $c=0$ cases, it was proved by Balas and Gauduchon \cite{BG} in 1985 for the Chern connection case, and by Sato and Sekigawa \cite{SS} in 1990 for the Riemannian connection case. The most difficult case is when $c>0$. This was resolved by Apostolov, Davidov and Muskarov \cite{ADM} in 1996, as a consequence of their classification on compact self-dual Hermitian surfaces. For dimensions $3$ or higher, the conjecture remain open in its full generality, with only partial results obtained in various special circumstances.

In \cite{DGM}, Davidov, Grantcharov, and Muskarov examined twistor spaces and showed that
\begin{theorem}[Davidov-Grantcharov-Muskarov\,\cite{DGM}]
Conjecture \ref{conj1} holds for all twistor spaces.
\end{theorem}

Note that if a twistor space has constant Chern (or Riemannian) holomorphic sectional curvature $c$,  then the constant $c$ has to be positive, and the twistor space has to be ${\mathbb C}{\mathbb P}^3$ equipped with (a constant multiple of) the standard metric. This gives a strong evidence towards the validity of the conjecture, as twistor spaces form a large class of compact complex threefolds.  

Another partial result was proved by K. Tang in \cite{Tang}, where he established Conjecture \ref{conj1} under the additional assumption that the metric is also Chern K\"ahler-like, meaning that the Chern curvature tensor $R^c$ obeys all K\"ahler symmetries, namely, $R^c_{i\overline{j}k\overline{\ell}} = R^c_{k\overline{j}i\overline{\ell}}$ for any $1\leq i,j,k,\ell \leq n$ where $n$ is the complex dimension of the manifold.  

\begin{theorem}[Tang\,\cite{Tang}]
Conjecture \ref{conj1} holds for all compact Chern K\"ahler-like manifolds.
\end{theorem}

In 2021, Chen, Chen, and Nie showed in \cite{CCN} that Conjecture \ref{conj1} is true if $c\leq 0$ and the manifold is locally conformally K\"ahler. Recall that a Hermitian manifold $(M^n,g)$ is said to be locally conformally K\"ahler (lcK for brevity) if $M$ is covered by a collection of open neighborhoods $\{ U_{\alpha}\}_{\alpha \in A}$, and in each  $U_{\alpha}$  there exists a smooth positive function $f_{\alpha}$ such that $f_{\alpha}g$ is a K\"ahler metric in $U_{\alpha}$. Note that given any  compact lcK manifold $(M^n,g)$, except the case when it is globally conformally K\"ahler (meaning that there exists a global positive smooth function $f$ such that $fg$ is K\"ahler), the manifold $M^n$ is always non-K\"ahlerian. The globally conformally K\"ahler case is said to be non-essential.  Essential lcK manifolds form an large and important subclass of non-K\"ahler manifolds. 

\begin{theorem}[Chen-Chen-Nie\,\cite{CCN}]
Conjecture \ref{conj1} holds for all compact locally conformally K\"ahler manifolds, provided that the constant $c$ in the statement of the conjecture is either negative or zero.
\end{theorem}

In the same paper, Chen, Chen, and Nie also constructed an explicit example which illustrated that the compactness assumption in Conjecture \ref{conj1} is a necessary one.

\begin{example}[{\bf 7}]
Let $g$ be the Hermitian metric on ${\mathbb C}^n$ with K\"ahler form
$$ \omega = \sqrt{-1} (1+|z|^2)^2 \partial \overline{\partial}\log  (1+|z|^2) = \sqrt{-1} \sum_{i,j=1}^n g_{i\bar{j}} \,dz_i \wedge d\overline{z}_j, \ \ \ \ g_{i\bar{j}}  = (1+|z|^2)\delta_{ij} - \overline{z}_iz_j,$$
where $(z_1, \ldots , z_n)$ is the standard coordinate of ${\mathbb C}^n$ and $|z|^2=\sum_{i=1}^n|z_i|^2$. The Chern curvature $R^c$ of $g$ is not identically zero, yet the Chern holomorphic sectional curvature $H^c=0$ identically, where
$ H^c(X)= R^c_{X\overline{X} X\overline{X}}/|X|^4  $
for any non-zero complex tangent vector $X$ of type $(1,0)$. Since $(g_{i\bar{j}})\geq I$, the metric $g$ dominates the complex Euclidean metric $g_0$, hence $({\mathbb C}^n, g)$ is complete. 
\end{example}

In \cite{ZhouZ}, the authors attempted to confirm Conjecture \ref{conj1} for all balanced threefolds. They were only able to prove the following weaker statement:

\begin{theorem}[Zhou-Zheng\,\cite{ZhouZ}] 
Any compact Hermitian threefold $(M^3,g)$ with vanishing real bisectional curvature must be Chern flat. 
\end{theorem}

Recall that {\em real bisectional curvature} is a curvature notion introduced by Yang and Zheng in \cite{XYangZ}. It is equivalent to $H^c$ in strength when the metric is K\"ahler, but it is slightly stronger than $H^c$ in the non-K\"ahler case. When the real bisectional curvature is identically zero, the metric is always  balanced. 

In \cite{RaoZ}, the authors attempted to examine Conjecture \ref{conj1} for a special class of pluriclosed metrics, namely, those which are Bismut K\"ahler-like, meaning that the curvature $R^b$ of the Bismut connection obeys all K\"ahler symmetries. See \cite{AOUV} for the exact definition of this interesting class. As a solution to the AOUV Conjecture, Zhao and Zheng proved in \cite{ZZ} that Bismut K\"ahler-like metrics are always pluriclosed.  

\begin{theorem}[Rao-Zheng\,\cite{RaoZ}] 
If a compact Hermitian manifold $(M^n,g)$ is Bismut K\"ahler-like and having constant Chern holomorphic sectional curvature, then $g$ is K\"ahler. 
\end{theorem}

In other words, the Chern connection case of Conjecture \ref{conj1} holds for compact Bismut K\"ahler-like manifolds. Of course it would be more desirable to confirm the conjecture for all pluriclosed manifolds, even just in dimension $3$.

A Hermitian manifold $(M^n,g)$ is called {\em Bismut torsion-parallel}, or BTP for brevity, if the torsion tensor $T^b$ of the Bismut connection $\nabla^b$ is parallel with respect to $\nabla^b$ itself, namely, if $\nabla^bT^b=0$. An equivalent form of AOUV conjecture, confirmed by Zhao and Zheng in \cite{ZZ},  states that Bismut K\"ahler-like metrics are always BTP. In fact BTP is actually a larger set than Bismut K\"ahler-like manifolds, as it contains all Vaisman manifolds. Recall that a Hermitian manifold is said to be {\em Vaisman} if it is locally conformally K\"ahler, and its Lee form is parallel under the Levi-Civita connection. When $n\geq 3$, Vasiman metrics and Bismut K\"ahler-like metrics form disjoint sets, both are subsets of BTP metrics. There are also balanced BTP metrics when $n\geq 3$.  We refer the readers to \cite{ZhaoZ24, ZhaoZ25} for more discussions on BTP manifolds. 

In \cite{ChenZ1}, Chen and Zheng discussed Conjecture \ref{conj1} on BTP metrics, and proved the following

\begin{theorem}[Chen-Zheng\,\cite{ChenZ1}] 
Let $(M^n,g)$ be a compact Hermitian manifold that is Bismut torsion-parallel (BTP). If either the metric $g$ is non-balanced or $n=3$, then Conjecture \ref{conj1} holds for $(M^n,g)$. 
\end{theorem}

As a corollary, we know that Conjecture \ref{conj1} holds for all compact Vaisman manifolds. Note that this corollary is contained in the theorem of Chen-Chen-Nie when $c\leq 0$, as Vaisman manifolds are special types of lcK manifolds. So the corollary is only new in the $c>0$ case. 

\vspace{0.2cm}

Next let us switch gear and consider Lie-Hermitian manifolds. A compact Hermitian manifold $(M^n,g)$ is said to be {\em Lie-Hermitian}, if $M=G/\Gamma$ where $G$ is a Lie group and $\Gamma \subset G$ a discrete subgroup, and the lift of the metric $g$ and the complex structure $J$ onto $G$ are both left-invariant. Note that we do not require $J$ to be bi-invariant, in other words $G$ may not be a complex Lie group. When $G$ is nilpotent or solvable, we will call $(M^n,g)$ a complex nilmanifold (or complex solvmanifold). Lie-Hermitian manifolds constitute a large and important class of compact locally homogeneous Hermitian manifolds. 

Note that the Lie group $G$ here must be unimodular as it admits a compact quotient. Denote by ${\mathfrak g}$ its Lie algebras. By a result of Hano \cite{Hano}, if a unimodular Lie algebra admits a K\"ahler structure, then it must be flat. So for Lie-Hermitian manifolds, the $c\neq 0$ cases in Conjecture \ref{conj1} do not occur, and the conjecture reduces to the following:

\begin{conjecture}
If the Chern (or Riemannian) holomorphic sectional curvature of a Lie-Hermitian manifold is a constant $c$, then $c=0$ and the metric must be Chern (or Riemannian) flat. 
\end{conjecture}

Riemannian flat Lie algebras were classified by Milnor in \cite{Milnor}. They are always unimodular and are abelian or $2$-step solvable. Such a Lie algebra, if even-dimensional, always admits a complex structure which makes it K\"ahler flat, and K\"ahler flat structures on Lie algebras  were explicitly characterized by Barberis, Dotti, and Fino in \cite{BDF} (see also \cite[Appendix]{VYZ} for a  slight variation). At this point of course one would naturally wonder if there exists any Riemannian flat Hermitian Lie algebras that is non-K\"ahler. To be more precise, let $({\mathfrak g},g)$ be an even-dimensional Lie algebra equipped with a flat metric. Since it is even-dimensional, there are lots of almost complex structures $J$ on the vector space ${\mathfrak g}$ which is compatible to $g$. The point here is that our $J$ is integrable, which means
$$ [x,y] - [Jx, Jy] + J[Jx, y] + J[x,Jy] =0, \ \ \ \ \ \ \forall \ x, y \in {\mathfrak g}. $$
This is equivalent to  the condition that the left-invariant almost complex structure on the corresponding Lie group is integrable. So the question now is if we have a Hermitian Lie algebra $({\mathfrak g},g,J)$ where $g$ is Riemannian flat, then must it be K\"ahler? This question was answered very recently by Zhang and Zheng:

 \begin{theorem}[Zhang-Zheng\,\cite{ZhangZ}] 
Any Riemannian flat Hermitian Lie algebra is K\"ahler flat. 
\end{theorem}

In contrast, Chern flat Lie algebras are still not well-understood yet. Complex Lie algebras (with any compatible metric) are certainly Chern flat, but there are also examples of Chern flat Lie algebras that are not complex Lie algebras. The classic theorem of Boothby \cite{Boothby} states that any compact Hermitian manifold that is Chern flat must be a quotient of a complex Lie group. However, it could be the quotient of another Lie group, namely, different (unisomorphic) Lie groups (with left-invariant complex structure and metric) could have their quotients be holomorphically isometric  to each other, and there are indeed such examples. So the classification of all unimodular Chern flat Hermitian Lie algebras is still a challenging task.

Coming back to Conjecture \ref{conj1} on Lie-Hermitian manifolds, the simplest situation is when the Lie group is nilpotent. The Chern case was proved by Li and Zheng in \cite{LiZ} for all complex nilmanifolds, while the Riemannian case was proved by Chen and Zheng in \cite[Theorem 2(1)]{ChenZ2} for complex nilmanifolds with nilpotent $J$ (in the sense of Cordero, Fern\'{a}ndez, Gray, and Ugarte \cite{CFGU}):

 \begin{theorem}[Li-Zheng\,\cite{LiZ}] 
If a complex nilmanifold has constant Chern holomorphic sectional curvature $c$, then $c=0$, it is Chern flat, and the Lie group is a complex Lie group.
\end{theorem}

 \begin{theorem}[Chen-Zheng\,\cite{ChenZ2}] 
If a complex nilmanifold has constant Riemannian holomorphic sectional curvature $c$ and $J$ is nilpotent in the sense of   \cite{CFGU}, then  $c=0$, the Lie group is abelian, and the manifold is  a flat complex torus. 
\end{theorem}

The success in confirming Conjecture \ref{conj1} for nilmanifolds are based on the important result of S. Salamon \cite{Salamon} which gives a nice basis for nilpotent Lie algebras with complex structure, along with a refined version by Cordero, Fern\'{a}ndez, Gray, and Ugarte in \cite{CFGU} when the complex structure $J$ is nilpotent in their sense. 

Currently  not much is known about Conjecture \ref{conj1} for Lie-Hermitian manifolds beyond the nilpotent case. It is not even known for all $2$-step solvable groups. Here we will mention two very special circumstances.

\begin{theorem}[Huang-Zheng\,\cite{HuangZ}] 
Let $(M^n,g)$ be a Lie-Hermitian manifold with universal cover $G$. Denote by ${\mathfrak g}$ the Lie algebra of $G$ and ${\mathfrak g}'=[{\mathfrak g},{\mathfrak g}]$ its commutator. Assume that ${\mathfrak g}$ is solvable and $J{\mathfrak g}'={\mathfrak g}'$. If the Chern holomorphic sectional curvature is constant, then it is Chern flat. 
\end{theorem}

\begin{theorem}[Li-Zheng\,\cite{LiZ2}] 
Let $(M^n,g)$ be a Lie-Hermitian manifold with universal cover $G$. Denote by ${\mathfrak g}$ the Lie algebra of $G$.   Assume that ${\mathfrak a}\subset {\mathfrak g}$ is an abelian ideal of codimension $2$ such that $J{\mathfrak a}={\mathfrak a}$. If $g$ has constant Chern (or Riemannian) holomorphic sectional  curvature, then it is Chern (or Riemannian) flat. 
\end{theorem}

In other words, Conjecture \ref{conj1} holds for all unimodular Hermitian Lie algebras which contain a $J$-invariant abelian ideal of codimension $2$. In particular, the conjecture holds for all {\em almost abelian} Lie-Hermitian manifolds. Recall that a Lie algebra is said to be almost abelian if it contains an abelian ideal of codimension $1$. In this case if ${\mathfrak b}$ is such an ideal, then $J{\mathfrak b}\cap {\mathfrak b}$ would be an abelian ideal of codimension $2$ which is $J$-invariant. 

\vspace{0.1cm}

Before we end the discussion on Conjecture \ref{conj1}, let us give a remark about this Hermitian space form problem for other canonical metric connections. Conjecture \ref{conj1} addressed the issue for the Chern and the Levi-Civita (Riemannian) connection. An obvious question is: what about the case for the Bismut connection $\nabla^b$? The analogous statement to Conjecture \ref{conj1} becomes the following pair:

\begin{conjecture} \label{conj6}
If a   compact Hermitian manifold $(M^n,g)$ has constant Bismut holomorphic sectional curvature $c$ and $c\neq 0$, then $g$ must be K\"ahler. 
\end{conjecture}

\begin{question} \label{question15}
Characterize all compact Hermitian manifolds which have vanishing Bismut holomorphic sectional curvature but are not Bismut flat. 
\end{question}

The reason for the diversion here is that there are actually examples of compact Hermitian manifolds with vanishing Bismut holomorphic sectional curvature but are not Bismut flat, e.g., the standard Hopf manifold $M^n=({\mathbb C}^n\setminus \{ 0\})/{\mathbb Z}f$ where $f(z)=\frac{1}{2}z$, equipped with the standard metric with K\"ahler form $\omega = \sqrt{-1}\, \partial \overline{\partial} |z|^2/|z|^2$. Its Bismut holomorphic sectional curvature is identically zero, but when $n\geq 3$ it is not Bismut flat. In \cite{ChenZ}, Chen and Zheng confirmed Conjecture \ref{conj6}  and answered Question \ref{question15} for the $n=2$ case. Their theorem states that, if $(M^2,g)$ is a compact Hermitian surface with constant Bismut holomorphic sectional curvature $c$, then $g$ is K\"ahler when $c\neq 0$, and when $c=0$ either $(M^2,g)$ is K\"ahler flat or $(M^2,g)$ is an isosceles Hopf surface equipped with an admissible metric. They give explicit description for admissible metrics. For higher dimensions, they also answered Conjecture \ref{conj6}  and  Question \ref{question15} for Bismut K\"ahler-like manifolds as well as complex nilmanifolds with nilpotent $J$. In \cite{ChenZ1}, Chen and Zheng extended this result to all BTP manifolds that are either non-balanced or of dimension $3$.

Besides the three canonical metric connections on Hermitian manifolds, namely the Chern, the Levi-Civita, and the Bismut connection,  there are also other metric connections that have been studied. For any real number $r$, the linear combination $D^r=\frac{1+r}{2}\nabla^c + \frac{1-r}{2}\nabla^b$ is called a {\em Gauduchon connection}, which is the line joining the Chern connection $\nabla^c=D^1$ and the Bismut connection $\nabla^b=D^{-1}$. More generally, one could consider the `plane' of canonical metric connections
\begin{equation}
D^r_s = (1-s)D^r+s\nabla, \ \ \ \ \ \ \ (r,s)\in \Omega := \{ s\neq 1\} \cup \{ (0,1)\} \subset {\mathbb R}^2, 
\end{equation}
where $\nabla$ is the Levi-Civita (Riemannian) connection. When $g$ is K\"ahler, all these connections coincide with $\nabla$, and when $g$ is not K\"ahler,  $D^r_s\neq D^{r'}_{s'}$ for any two different points $(r,s)$ and $(r',s')$ in $\Omega$. In the spirit of Conjecture \ref{conj1}, a natural question to ask is the following: 

\begin{question}  \label{question16}
Let $(M^n,g)$ be a compact Hermitian manifold and fix a point  $(r,s)\in \Omega$.

(1). When will it be $D^r_s$-flat, or more generally when will it be $D^r_s$-K\"ahler-like (namely the curvature of $D^r_s$ obeys all K\"ahler symmetries)?

(2). When will the holomorphic sectional curvature of $D^r_s$ be a constant $c$? 
\end{question}

For part (1), Lafuente and Stanfield \cite{LS} give the following beautiful characterization for all Gauduchon connections (other than Chern or Bismut)
\begin{theorem}[Lafuente-Stanfield\,\cite{LS}]
For any $r\neq \pm 1$, if the Gauduchon connection $D^r$ of  a compact Hermitian manifold  $(M^n,g)$ is K\"ahler-like, then $g$ must be K\"ahler.
\end{theorem}
Prior to that, partial results were obtained by Yang-Zheng \cite{YZ1} and Fu-Zhou \cite{FuZhou}, where a critical interval of $r$-values was excluded. Away from the Gauduchon line, the question was answered by Zhao and Zheng in \cite{ZhaoZ}. To state their result, let
$$ \Omega'=\Omega \setminus \big(  \{ s=0\} \cup \{ (1,0), (0,-1), (-1,2), (\frac{1}{3},-2)\} \big)  $$
be the region in the parameter space excluding the Gauduchon line along with four special points. The corresponding $D^r_s$ for these four points are denoted respectively by $\nabla$, $\nabla'$, $\nabla^+$, and $\nabla^-$. Note that $\nabla$ is the Levi-Civita connection, and $\nabla'$ is the reflection point of $\nabla$ with respect to the Gauduchon line, while $\nabla^+$, $\nabla^-$ are two special connections closely related to the Bismut connection $\nabla^b$. Their result states that

\begin{theorem}[Zhao-Zheng\,\cite{ZhaoZ}]
For any $(r,s)\in \Omega'$, if the canonical metric connection $D^r_s$ of a compact Hermitian manifold  $(M^n,g)$ is K\"ahler-like, then $g$ must be K\"ahler. Furthermore, $\nabla'$ will be K\"ahler-like if and only if $\nabla$ is K\"ahler-like, while $\nabla^+$ or $\nabla^-$ will be K\"ahler-like if and only if $\nabla^b$ is K\"ahler-like.
\end{theorem}

Now let us turn our attention to part (2) of Question \ref{question16}. In \cite{CN}, Chen and Nie discovered a particular curve in the parameter space which plays a crucial role in the space form problem:
\begin{equation}
\Gamma = \{ (r,s)\in {\mathbb R}^2  \mid (1-r+rs)^2 + s^2=4 \} \subset \Omega .
\end{equation}
We will call $\Gamma$ the {\em Chen-Nie curve} from now on. The intersection  of $\Gamma$ with the Gauduchon line $\{ s=0\}$ is the point  $(-1,0)$, which corresponds to the Bismut connection. Note that $\Gamma$ also passes through both $\nabla^+$ and $\nabla^-$. Chen and Nie formulated the following:

\begin{conjecture}[Chen-Nie\,\cite{CN}]
Let $(M^n,g)$ be a compact Hermitian manifold. Assume that the holomorphic section curvature of $D^r_s$ is a constant $c$. If $c\neq 0$, then $g$ must be K\"ahler. If $c=0$ and $(r,s)\in \Omega \setminus \Gamma$, then $g$ must be $D^r_s$-flat.  
\end{conjecture}

Note that the in last sentence, the $D^r_s$-flatness would imply that $g$ is K\"ahler for all $(r,s)\in \Omega \setminus \Gamma$ except $(0,1)$ and $(0,-1)$, by the result of Zhao and Zheng mentioned above. For points on $\Gamma$, one could also raise the following

\begin{question}
For $(r,s)\in \Gamma$, what kind of compact Hermitian manifolds will have its $D^r_s$ connection being non-flat but with vanishing holomorphic sectional curvature?
\end{question}

Chen and Nie in \cite{CN} completely resolved the $2$-dimensional  case:

\begin{theorem}[Chen-Nie\,\cite{CN}]
Let $(M^2,g)$ be a compact Hermitian surface with pointwise constant
 holomorphic sectional curvature with respect to its $D^r_s$ connection. Then either $g$ must be K\"ahler,
 or $(r,s) \in \Gamma$ and $(M^2,g)$ is an isosceles Hopf surface equipped with an admissible metric.
\end{theorem}

For higher dimensions, Chen and Zheng in \cite{ChenZ2} give some partial answers to the above conjecture of Chen-Nie in several special situations:

\begin{theorem}[Chen-Zheng\,\cite{ChenZ2}]
Let $(M^n,g)$ be a compact Hermitian manifold. Assume that the holomorphic sectional curvature of $D^r_s$ is a constant $c$. 

(1). If $(M^n,g)$ is a complex nilmanifold with nilpotent $J$ in the sense of \cite{CFGU}, then $c=0$. Furthermore, if $D^r_s\neq \nabla^c$ then the Lie group $G$ is abelian and $g$ is K\"ahler flat; if $D^r_s= \nabla^c$ then $G$ is a complex Lie group and $g$ is Chern flat.

(2). If $(M^n,g)$ is a non-balanced Bismut torsion-parallel manifold, then $c=0$ and $(r,s)\in \Gamma$.

(3). If  $(M^3,g)$ is a balanced Bismut torsion-parallel manifold, then either $g$ is K\"ahler, or $c=0$ and $D^r_s=\nabla^c$ and $g$ is Chern flat. 
\end{theorem}

In other words, the Chen-Nie Conjecture  holds for all complex nilmanifolds with nilpotent $J$, all compact non-balanced Bismut torsion-parallel manifolds, and all compact balanced Bismut torsion-parallel threefolds. 

\vspace{0.3cm}

\section{Fino-Vezzoni conjecture}

In this section, we will discuss the history and recent development towards the Fino-Vezzoni conjecture,  stated as Conjecture \ref{conj2} in the introduction section. The conjecture says that any compact non-K\"ahlerian manifold $M^n$ cannot simultaneously admit a balanced metric and a pluriclosed metric. 

First let us recall the definitions. Given a compact Hermitian manifold $(M^n,g)$, denote by $\omega$ its K\"ahler form. Then $g$ is {\em balanced} if $\partial \overline{\partial} (\omega^{n-1})=0$, and $g$ is {\em pluriclosed} if $\partial \overline{\partial} \omega =0$. Since K\"ahlerness means $d\omega =0$, both the balanced condition and the pluriclosed condition are generalizations to the K\"ahlerness. A few quick remarks. First, a theorem of Alexsandrov and Ivanov \cite{AI} says that if a metric $g$ is balanced and pluriclosed at the same time, then it must be K\"ahler. So in the statement of Conjecture \ref{conj2} it is vital to allow the two special Hermitian metrics to be different.  Secondly, since the balanced condition coincides with K\"ahlerness when $n=2$, so the starting dimension of Conjecture \ref{conj2} is $n=3$. Unfortunately, at present the conjecture is still open in its full generality even for $n=3$. Thirdly, the compactness assumption in Conjecture \ref{conj2} is generally believed to be a necessary one, namely, Conjecture \ref{conj2} might fail if $M^n$ is non-compact, though such a counterexample could be difficult to construct. 

To illustrate the subtlety of this issue, let us consider Lie groups equipped with left-invariant complex structures. In this case we only consider left-invariant metrics. In \cite{FSwann}, Freibert and Swann constructed explicit examples of a $2$-step solvable Lie group $G$ equipped with a left-invariant complex structure $J$, such that $(G,J)$ admits compatible left-invariant metrics $g$ and $h$ so that $g$ is balanced and $h$ is pluriclosed, yet $(G,J)$ does not admit any compatible left-invariant metric that is K\"ahler. See Example 4.4 and Example 4.16 in \cite{FSwann}. Of course the Lie groups $G$ involved here are not {\em unimodular,} hence they cannot admit any compact quotient. 

Balanced metrics and pluriclosed metrics (the latter is also called SKT metrics in some literature, which stands for {\em strong K\"ahler with torsion}) are both extensively studied in complex geometry. For instance, the balanced structure on compact complex manifolds are known to be preserved under bimeromorphism:

\begin{theorem}[Alessandrini-Bassanelli\,\cite{AB}]
Let $X$ and $Y$ be compact complex manifolds that are bimeromorphic to each other. If $X$ admits a balanced metric, then so does $Y$. 
\end{theorem}

So for bimeromorphic changes the balanced condition is more robust than K\"ahlerness. On the other hand, the balanced condition is not preserved under small deformation, so in that sense it is not as `stable' as K\"ahlerness. It is still an active research topic to find out assumptions (e.g., various forms of $\partial\overline{\partial}$-lemma) under which the balanced condition is preserved under deformation. 

An important source of interest on balanced metrics comes from string theory. The Hull-Strominger system depicts a model for the hidden space as a compact complex threefold which admits a metric that is conformal to a (non-K\"ahler) balanced metric, plus other restrictions. The fact that the naturally induced metrics $h_t$ on all twistor spaces are always balanced indicates that there are lots of non-K\"ahler balanced manifolds. We refer the readers to \cite{Fu} and the references therein for more discussion on balanced metrics. 

Pluriclosed metrics are also widely studied in the past a couple of decades. They appear naturally in the Hermitian curvature flow theory of Streets and Tian, see \cite{StreetsTian, StreetsTian2} for the discussion on Hermitian curvature flow and in particular the pluriclosed flow. There is an excellently written survey paper by Fino and Tomassini \cite{FinoTomassini} which gives a thorough  discussion on this special type of Hermitian structure. 

Conjecture \ref{conj2} predicts that in the non-K\"ahler universe, balanced structure and pluriclosed structure cannot co-exist on the same manifold. So far the conjecture is still wide open in its full generality, even in the lowest dimension $n=3$ case. However, there are lots of partial results which confirm the conjecture in many special situations, giving strong evidence towards the validity of the conjecture. For instance, Fu, Li, and Yau in \cite{FuLiYau} proved that the conjecture holds for a special type of non-K\"ahler Calabi-Yau threefolds, and Fei \cite{Fei} confirmed the conjecture for generalized Calabi-Gray manifolds. We will omit the detailed statements here to avoid technicalities, and refer the readers to the references given above.  

\begin{theorem}[Verbitsky\,\cite{Verbitsky}]
All non-K\"ahlerian twistor space admits no pluriclosed metric. 
\end{theorem}

In particular,  Conjecture \ref{conj2} holds for all twistor spaces. This provides a strong piece of evidence towards the conjecture.

\begin{theorem}[Chiose\,\cite{Chiose}]
If $M^n$ is bimeromorphic to a compact K\"ahler manifold, and $M^n$ admits a pluriclosed metric, then it will admit a K\"ahler metric.
\end{theorem}

The set of all compact complex manifolds bimeromorphic to compact K\"ahler manifolds are called the {\em Fujiki C-class}. The above theorem of Chiose says that Conjecture \ref{conj2} holds for the Fujiki C-class. 

The Oeljeklaus-Toma manifolds is a special class of compact non-K\"ahler manifolds which serves as the high dimensional analogue of a type of Inoue surfaces. In a more recent work \cite{Otiman}, Otiman proved that the conjecture holds for all Oeljeklaus-Toma manifolds:

\begin{theorem}[Otiman\,\cite{Otiman}]
Conjecture \ref{conj2} holds for all Oeljeklaus-Toma manifolds.
\end{theorem}

\begin{theorem}[Angella-Otiman\,\cite{AngellaOtiman}]
Any compact Vaisman manifold of complex dimension at least $3$ does not admit any balanced or pluriclosed metric. In particular, Conjecture \ref{conj2} holds for Vaisman manifolds.
\end{theorem}

\begin{theorem}[Zhao-Zheng\,\cite{ZZ}]
Any compact, non-K\"ahler, Bismut K\"ahler-like manifold does not admit any balanced metric. In particular, Conjecture \ref{conj2} holds for all compact Bismut K\"ahler-like  manifolds.
\end{theorem}

What the authors proved is a slightly stronger result, namely, if $(M^n,g)$ is a non-K\"ahler, compact Bismut K\"ahler-like manifold, then $M^n$ does not admit any {\em strongly Gauduchon metric}, which means a Hermitian metric $h$ on $M^n$ whose K\"ahler form $\omega_h$ satisfies the condition that $\partial (\omega_h^{n-1})=\overline{\partial}\Phi$ for some global $(n,n-2)$-form $\Phi$ on $M^n$. Since balanced metrics are certainly strongly Gauduchon, we know that the manifold $M^n$ in the above theorem cannot admit any balanced metric.

Next let us consider Lie-Hermitian manifolds. The nilpotent case was confirmed by Fino and Vezzoni \cite{FV15,FV16}: 

\begin{theorem}[Fino-Vezzoni\,\cite{FV15, FV16}]
Conjecture \ref{conj2} holds for all complex nilmanifolds, and it also holds for all real $6$-dimensional complex solvmanifolds of Calabi-Yau type.
\end{theorem}

For the nilmanifold case, they implicitly assumed that the nilpotent group is of step at most $2$. Later on,  Arroyo and Nicolini showed in \cite{ArroyoN} that this assumption can be removed, as for any complex nilmanifold admitting a pluriclosed metric, the corresponding Lie group has to be $2$-step nilpotent (or abelian). 

In the following, we will always assume that $(M^n,g)$ is a Lie-Hermitian manifold, namely, a compact Hermitian manifold with $M=G/\Gamma$, where $G$ is a Lie group, $\Gamma \subset G$ is a discrete subgroup, and both the metric $g$ and the complex structure $J$ (when lifted onto $G$) are left-invariant. Denote by ${\mathfrak g}$ the Lie algebra of $G$. Beyond nilmanifolds, Conjecture \ref{conj2} was confirmed in the following special cases:

\begin{theorem}[Fino-Grantcharov-Vezzoni\,\cite{FGV}, Podest\`a\,\cite{Podesta}]
Conjecture \ref{conj2} holds for all compact semi-simple Lie groups. 
\end{theorem}

\begin{theorem}[Giusti-Podest\`a\,\cite{GPodesta}]
Conjecture \ref{conj2} holds for all Lie-Hermitian manifolds where $G$ is non-compact semi-simple and $J$ is regular.  
\end{theorem}

\begin{theorem}[Fino-Paradiso\,\cite{FP,FP1,FP2}]
Conjecture \ref{conj2} holds for all Lie-Hermitian manifolds where $G$ is almost abelian or almost nilpotent solvable.
\end{theorem}

Recall that a Lie algebra ${\mathfrak g}$  is said to be {\em almost abelian} if it contains an abelian ideal of codimension $1$. It is said to be {\em almost nilpotent solvable}, if it is solvable and its nilradical has $1$-dimensional commutator. A simple and important class beyond the nilpotent ones are $2$-step solvable Lie algebras, namely, a non-abelian Lie algebra ${\mathfrak g}$ whose commutator ${\mathfrak g}'=[{\mathfrak g}, {\mathfrak g}]$  is abelian. It would certainly be highly desirable to confirm Conjecture \ref{conj2} for this class. In their recent work \cite{FSwann22, FSwann}, Freibert and Swann systematically studied the Hermitian structures on $2$-step solvable Lie algebras, and they proved  Conjecture \ref{conj2} for a large subset, namely when the complex structure $J$ is in one of the following three `pure types':
\begin{itemize}
\item Pure type I: $J{\mathfrak g}'\cap {\mathfrak g}'=0$;
\item Pure type II: \,$J{\mathfrak g}' = {\mathfrak g}'$;
\item Pure type III: $J{\mathfrak g}'+ {\mathfrak g}'={\mathfrak g}$.
\end{itemize} 

\begin{theorem}[Freibert-Swann\,\cite{FSwann22, FSwann}]
Conjecture \ref{conj2} holds for all Lie-Hermitian manifolds where $G$ is  $2$-step solvable and $J$ is of pure types.
\end{theorem}

Now let us consider a natural generalization to the almost abelian case, namely when the Lie algebra ${\mathfrak g}$ contains an abelian ideal ${\mathfrak a}$ of codimension $2$. Suppose $(J,g)$ is a Hermitian structure on ${\mathfrak g}$. Then ${\mathfrak a}_J = {\mathfrak a} \cap J{\mathfrak a}$ is a $J$-invariant  ideal of ${\mathfrak g}$ of codimension either $2$ or $4$, depending on whether $J{\mathfrak a}={\mathfrak a}$ or not. Note that ${\mathfrak g}$ is always solvable of step at most $3$, but it will not be of step $2$ in general. For the easier $J{\mathfrak a}={\mathfrak a}$ case, Fino-Vezzoni Conjecture was verified by Li and Zheng in \cite{LiZ1}:

\begin{theorem}[Li-Zheng\,\cite{LiZ1}]
Conjecture \ref{conj2} holds for all Lie-Hermitian manifolds if the Lie algebra ${\mathfrak g}$ of $G$ contains a $J$-invariant abelian ideal of codimension $2$.
\end{theorem}

Last year, utilizing the earlier work by Guo and Zheng in \cite{GuoZ}, Cao and Zheng in \cite{CaoZ2} extended the above theorem by removing the $J{\mathfrak a}={\mathfrak a}$ assumption:

\begin{theorem}[Cao-Zheng\,\cite{CaoZ2}]
Conjecture \ref{conj2} holds for all Lie-Hermitian manifolds if the Lie algebra ${\mathfrak g}$ of $G$ contains an abelian ideal of codimension $2$.
\end{theorem}

\vspace{0.3cm}

\noindent\textbf{Acknowledgments.} We are very grateful to Professor Yau for his constant support and encouragement throughout all these years. We would like to thank mathematicians Huai-Dong Cao, Haojie Chen, Lei Ni, Xiaolan Nie, Kai Tang, Bo Yang, Xiaokui Yang, and Quanting Zhao for their interest and helpful discussions. We are also thankful to our graduate students Kexiang Cao, Shuwen Chen, Yuqin Guo, Xin Huang, Yulu Li, Peipei Rao, and Dongmei Zhang for the fruitful collaborations, which are well-reflected in this survey. Finally, we would like to express our sincere gratitude to all three referees for their valuable suggestions and corrections, which improved the readability of the article. 

\vs

\end{document}